\documentclass[graybox]{svmult}

\usepackage{type1cm}        
\usepackage{makeidx}         
\usepackage{graphicx}        
\usepackage{multicol}        
\usepackage[bottom]{footmisc}

\usepackage{newtxtext}       %
\usepackage[varvw]{newtxmath}       

\usepackage[numbers,sort&compress]{natbib}
\usepackage{float}

\usepackage{xcolor}
\usepackage{scalerel}

\makeindex             

\begin{document}

\title*{Optimized Schwarz Waveform Relaxation\\ for the Damped Wave Equation}
\author{Gerardo Cicalese\orcidID{0000-0001-5061-8608} and\\ 
Gabriele Ciaramella\orcidID{0000-0002-5877-4426} and\\ 
Ilario Mazzieri\orcidID{0000-0003-4121-8092} and\\ 
Martin J.~Gander\orcidID{0000-0001-8450-9223}}
\institute{Gerardo Cicalese, Gabriele Ciaramella, Ilario Mazzieri \at MOX Laboratory, Dipartimento di Matematica, Politecnico di Milano, Piazza Leonardo da Vinci 32, 20133 Milano, Italia, \email{{gerardo.cicalese, gabriele.ciaramella, ilario.mazzieri}@polimi.it}
\and Martin J.~Gander \at Section de mathématiques, Université de Genève, Rue du Conseil-Général 9, CP 64, 1211 Genève 4, Suisse, \email{martin.gander@unige.ch}}
%
%
\maketitle

\abstract*{The performance of Schwarz Waveform Relaxation is critically dependent on the choice of transmission conditions. While classical absorbing conditions work well for wave propagation, they prove insufficient for damped wave equations, particularly in viscoelastic damping regimes where convergence becomes prohibitively slow. This paper addresses this limitation by introducing a more general transmission operator with two free parameters for the one-dimensional damped wave equation. Through frequency-domain analysis, we derive an explicit expression for the convergence factor governing the convergence rate. We propose and compare two optimization strategies---\texorpdfstring{$L_\infty$}{L-infinity} and \texorpdfstring{$L_2$}{L-2} minimization---for determining optimal transmission parameters. Numerical experiments demonstrate that our optimized approach significantly accelerates convergence compared to standard absorbing conditions, especially for viscoelastic damping cases. The method provides a computationally efficient alternative to exhaustive parameter search while maintaining robust performance across different damping regimes.}

\section{Introduction}
The development of transmission conditions in overlapping Schwarz Waveform Relaxation (SWR) methods for the undamped wave equation has progressed from classical Dirichlet to absorbing conditions. While Dirichlet conditions produce strong interface reflections that slow down convergence, absorbing conditions  reduce error reflections and promote efficient wave transmission between subdomains \cite{engquist_absorbing_1977,gander_optimal_1999,gander_absorbing_2004}.

In this paper, we investigate general transmission conditions for the (one-dimensional) damped wave equation, for which no practical transparent conditions are available in the literature.
In particular, for a domain $\Omega = (0,L)$ and time $t \geq0$, we consider the viscoelastic-telegrapher damped wave equation
\begin{equation} \label{eq:1d_wave_eq}
    \partial_{tt} u + \gamma \partial_t u = c^2 \partial_{xx} u + \nu \partial_t \partial_{xx} u + f, \quad \text{in} \; \Omega,
\end{equation}
with initial conditions $u(x,0) = u_0(x), \; \partial_t u(x,0) = v_0(x), \; x \in \Omega$, and homogeneous Dirichlet boundary conditions $u(0,t) = u(L,t) = 0, \; t > 0$.
Here, $\gamma$ represents the telegrapher damping and $\nu$ the viscoelastic damping. Problem \eqref{eq:1d_wave_eq} is hyperbolic when $\nu = 0$ \cite{zauderer_partial_2006} and not hyperbolic when $\nu > 0$~\cite{liu_spectrum_1998}.

For this model, we investigate transmission conditions characterized by two parameters that can be tuned to optimize SWR convergence. We consider a spectral optimization framework that automatically identifies the optimal parameters. This relies on a frequency-domain analysis to derive a closed-form expression for the SWR convergence factor, which we then minimize (numerically) using two different metrics. The proposed approach provides a systematic way to tune the transmission conditions. It turns out that our approach significantly accelerates the SWR convergence in the case of viscoelastic damping.

To define the SWR method for our model problem \eqref{eq:1d_wave_eq}, we decompose $\Omega = (0,L)$ into two overlapping subdomains $\Omega_1 = (0,b)$ and $\Omega_2 = (a,L)$ with $0 < a < b < L$ and overlap $\delta = b - a > 0$. The SWR iterations are then defined as
\begin{align}
    \partial_{tt} v^{(k+1)} + \gamma \partial_t v^{(k+1)} &= c^2 \partial_{xx} v^{(k+1)} + \nu \partial_t \partial_{xx} v^{(k+1)} + f, & \text{in} \; \Omega_1, \\
    \big(\partial_x + \Lambda\big) v^{(k+1)}(b,t) &= \big(\partial_x + \Lambda\big) w^{(k)}(b,t), & \label{eq:interface_conditions_v} \\
    \partial_{tt} w^{(k+1)} + \gamma \partial_t w &= c^2 \partial_{xx} w^{(k+1)} + \nu \partial_t \partial_{xx} w^{(k+1)} + f, & \text{in} \; \Omega_2, \\
    \big(\partial_x - \Lambda\big) w^{(k+1)}(a,t) &= \big(\partial_x - \Lambda\big) v^{(k)}(a,t), \label{eq:interface_conditions_w}
\end{align}
where $k$ is the iteration index and $\Lambda = p \partial_t + q$ with $p, q \in \mathbb{R}$. Notice that the classical transparent condition for $\gamma = \nu = 0$ corresponds to $(p, q) = (1/c, 0)$ \cite{engquist_absorbing_1977,gander_optimal_1999}.

The paper is organized as follows. In Section \ref{sec:2}, we derive the convergence factor and introduce two optimization strategies. In Section \ref{sec:3}, we present numerical experiments, including convergence analyses for telegrapher and viscoelastic damping with various parameters, and we examine how optimized transmission-condition parameters vary with the damping coefficients. We present our conclusions in Section \ref{sec:conclusions}.

\section{Convergence Factor and Optimization Strategies} \label{sec:2}
For the error analysis, we set $f \equiv 0$, use the Laplace transform in time, and denote by $\widehat{u}(x,s)$ the Laplace-transformed field. From \eqref{eq:1d_wave_eq} we get 
\begin{equation}\label{eq:Laplace}
    \partial_{xx} \widehat{u} = - \kappa^2(s) \widehat{u},
\end{equation}
where $\kappa(s) := -\mathrm{i} \, s/c \sqrt{ \left( 1 + \gamma s^{-1} \right) / \left( 1 + \nu c^{-2} s \right)}$ and $\mathrm{i}$ is the imaginary unit.
To derive the convergence factor, we first solve \eqref{eq:Laplace} in each subdomain and apply the homogeneous Dirichlet boundary conditions to get
\begin{align}
    \widehat{v}(x,s) &= A_v(s) \left(e^{\mathrm{i} \kappa(s) x} - e^{-\mathrm{i} \kappa(s) x}\right), \label{eq:traveling_wave_sols_v} \\
    \widehat{w}(x,s) &= A_w(s) \left(e^{\mathrm{i} \kappa(s) x} - e^{2\mathrm{i}\kappa L} e^{-\mathrm{i} \kappa(s) x}\right). \label{eq:traveling_wave_sols_w}
\end{align}
By transforming the transmission conditions \eqref{eq:interface_conditions_v} and \eqref{eq:interface_conditions_w}, we get
\begin{align}
    (\partial_x + \lambda(s)) \; \widehat{v}^{(k+1)}(b,s) &= (\partial_x + \lambda(s)) \; \widehat{w}^{(k)}(b,s), \label{eq:interface_conditions_laplace_v} \\
    (\partial_x - \lambda(s)) \; \widehat{w}^{(k+1)}(a,s) &= (\partial_x - \lambda(s)) \; \widehat{v}^{(k)}(a,s), \label{eq:interface_conditions_laplace_w}
\end{align}
where $\lambda(s) = p s + q$.
Inserting \eqref{eq:traveling_wave_sols_v} and \eqref{eq:traveling_wave_sols_w} into \eqref{eq:interface_conditions_laplace_v} and \eqref{eq:interface_conditions_laplace_w}, we obtain
\begin{align}
    A_v^{(k+1)}(s) D_v(s) &= A_w^{(k)}(s) N_w(s), \label{eq:if_b_amp} \\
    A_w^{(k+1)}(s) D_w(s) &= A_v^{(k)}(s) N_v(s), \label{eq:if_a_amp}
\end{align}
where
\begin{align*}
    D_v(s) &= \mathrm{i}\kappa(s) \left(e^{\mathrm{i} \kappa(s) b} + e^{-\mathrm{i} \kappa(s) b}\right)
    + \lambda(s) \left(e^{\mathrm{i} \kappa(s) b} - e^{-\mathrm{i} \kappa(s) b}\right),
    \\ N_w(s) &= \mathrm{i}\kappa(s) \left(e^{\mathrm{i} \kappa(s) b} + e^{2\mathrm{i} \kappa(s) L}e^{-\mathrm{i} \kappa(s) b}\right)
    + \lambda(s) \left(e^{\mathrm{i} \kappa(s) b} - e^{2\mathrm{i} \kappa(s) L}e^{-\mathrm{i} \kappa(s) b}\right),
    \\ D_w(s) &= \mathrm{i} \kappa(s) \left(e^{\mathrm{i} \kappa(s) a} + e^{2\mathrm{i} \kappa(s) L}e^{-\mathrm{i} \kappa(s) a}\right)
    - \lambda(s) \left(e^{\mathrm{i} \kappa(s) a} - e^{2\mathrm{i} \kappa(s) L}e^{-\mathrm{i} \kappa(s) a}\right),
    \\ N_v(s) &= \mathrm{i} \kappa(s) \left(e^{\mathrm{i} \kappa(s) a} + e^{-\mathrm{i} \kappa(s) a}\right)
    - \lambda(s) \left(e^{\mathrm{i} \kappa(s) a} - e^{-\mathrm{i} \kappa(s) a}\right).
\end{align*}
Combining \eqref{eq:if_b_amp} and \eqref{eq:if_a_amp} yields the two-step recurrence
\begin{equation*}
    A_v^{(k+1)}(s) = G^2(s;p,q) A_v^{(k-1)}(s), \; \text{ where } \; G^2(s;p,q) = \frac{N_w(s) \, N_v(s)}{D_v(s) \, D_w(s)}.
\end{equation*}
The convergence factor is therefore given after simplification by
\begin{equation*}
    \rho(s; p, q) = |G(s;p,q)| = \left| \frac{\kappa(s) \, \cos(\kappa(s) \, a) - (p s + q) \, \sin(\kappa(s) \, a)}{\kappa(s) \, \cos(\kappa(s) \, b) + (p s + q) \, \sin(\kappa(s) \, b)} \right|.
\end{equation*}
Next, we consider two optimization strategies for the parameters $(p, q)$, namely, spectral $L_\infty$ and $L_2$ optimization. The optimization is performed on the imaginary axis $s=\mathrm{i} \omega$, since the relevant modes of the wave equation are oscillatory. The goal is to find parameters minimizing the values of $\rho$ across the relevant frequency band $[\omega_{\min},\omega_{\max}]$. Since the total simulation time is $T$, the smallest resolvable nonzero frequency is $\omega_{\min} = \pi/T$; see \cite[Figure 3.17]{gander_time_2024}. Due to the discrete time step $\Delta t$, the largest frequency that can be represented without aliasing is the Nyquist frequency $\omega_{\max} = \pi/\Delta t$. For the $L_\infty$ optimization, we minimize the global convergence factor $\hat{\rho}_{\infty}$, i.e.,
\begin{equation*}
    \min_{(p, q) \in \mathbb{R}^2} \hat{\rho}_{\infty}, \quad \text{where} \quad \hat{\rho}_{\infty} =
    \max_{\omega \in [\omega_{\min}, \omega_{\max}]} 
    \rho(\mathrm{i} \omega; p, q),
\end{equation*}
while for the $L_2$ optimization, we minimize the global convergence factor $\hat{\rho}_2$, i.e.,
\begin{equation*}
    \min_{(p, q) \in \mathbb{R}^2} \hat{\rho}_2, \quad \text{where} \quad \hat{\rho}_2 =
    \sqrt{\frac{1}{\omega_{\max}-\omega_{\min}}
    \int_{\omega_{\min}}^{\omega_{\max}} \rho(\mathrm{i} \omega;p,q)^2 d\omega}.
\end{equation*}
Note that $\hat{\rho}_2$ represents a normalized $L_2$ norm, in particular the root-mean-square (RMS) value; this normalization ensures that $\hat{\rho}_2$ remains on the same scale as $\hat{\rho}_{\infty}$, without altering the location of the minimizer.
To perform the optimization numerically, we discretize the frequency interval $[\omega_{\min}, \omega_{\max}]$ using a grid of 1000 nodes. For the $L_2$ criterion, the integral is approximated using the trapezoidal rule (implemented as \texttt{trapz} in MATLAB). The minimization is carried out using the Nelder-Mead simplex method (\texttt{fminsearch}), with termination tolerances set to $10^{-4}$ for both the parameter values (\texttt{TolX}) and the objective function (\texttt{TolFun}).

\section{Numerical Experiments} \label{sec:3}
To validate the parameters obtained by the optimization frameworks discussed in Section \ref{sec:2}, we compare them to those minimizing the SWR error after $k=10$ iterations. To do so, we discretize \eqref{eq:1d_wave_eq} by the Finite Difference Time-Domain (FDTD) scheme
\vspace{-0.3cm}
\begin{multline*}
    \frac{u_i^{\,n+1}-2u_i^{\,n}+u_i^{\,n-1}}{\Delta t^{2}}
    +\frac{\gamma}{2\Delta t}\bigl(u_i^{\,n+1}-u_i^{\,n-1}\bigr)
    =
    \frac{c^{2}}{\Delta x^{2}}\bigl(u_{i+1}^{\,n}-2u_i^{\,n}+u_{i-1}^{\,n}\bigr)
    \\
    +
    \frac{\nu}{2\Delta t\,\Delta x^{2}}
    \Bigl[(u_{i+1}^{\,n+1}-2u_i^{\,n+1}+u_{i-1}^{\,n+1})
    -(u_{i+1}^{\,n-1}-2u_i^{\,n-1}+u_{i-1}^{\,n-1})\Bigr]
    + f_i^{\,n}.
\end{multline*}
Notice that, for $\nu = 0$, the scheme becomes explicit. To discretize the transmission conditions, we use second-order stencils for both spatial (at $x = b$ and $x = a$) and time derivatives. We choose the following parameters: wave speed $c = 1$, subdomain interface position $a = 0.3$, overlap length $\delta = 0.1$, spatial discretization $\Delta x = 0.002$, temporal discretization $\Delta t = 0.002$, and final time $T = 5$. For these parameters, the discretization error is less than $0.1\%$. The error between the computed solution $\mathbf{u}_{\text{SWR}}$ and the reference solution $\mathbf{u}_{\text{ref}}$ (FDTD discrete solution) is computed using the relative infinity norm at $t = T$: $\|\mathbf{u}_{\text{SWR}} - \mathbf{u}_{\text{ref}}\|_\infty / \|\mathbf{u}_{\text{ref}}\|_\infty$.

\begin{figure}[t]
    \centering

    \centering
    \begin{minipage}{\linewidth}
        \centering
        \includegraphics[width=1\linewidth,
        trim=4.6cm 1.2cm 3.8cm 0.5cm, clip]{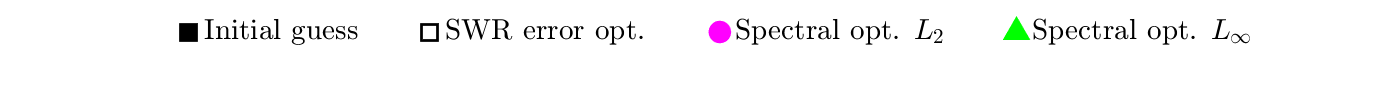}
    \end{minipage}

    \begin{minipage}[c]{0.48\linewidth}
        \centering
        \includegraphics[width=\linewidth,
        trim=0.3cm 0cm 1.3cm 0.8cm, clip]{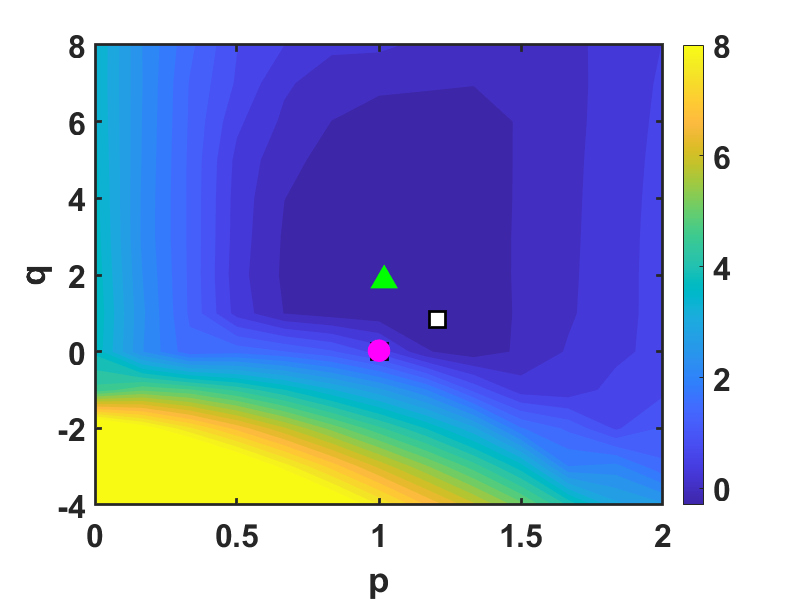}
        \textbf{(a$_\gamma$: $\gamma = 4, \nu = 0$)}
    \end{minipage}%
    \hfill
    \begin{minipage}[c]{0.48\linewidth}
        \centering
        \includegraphics[width=\linewidth,
        trim=0.3cm 0cm 1.3cm 0.8cm, clip]{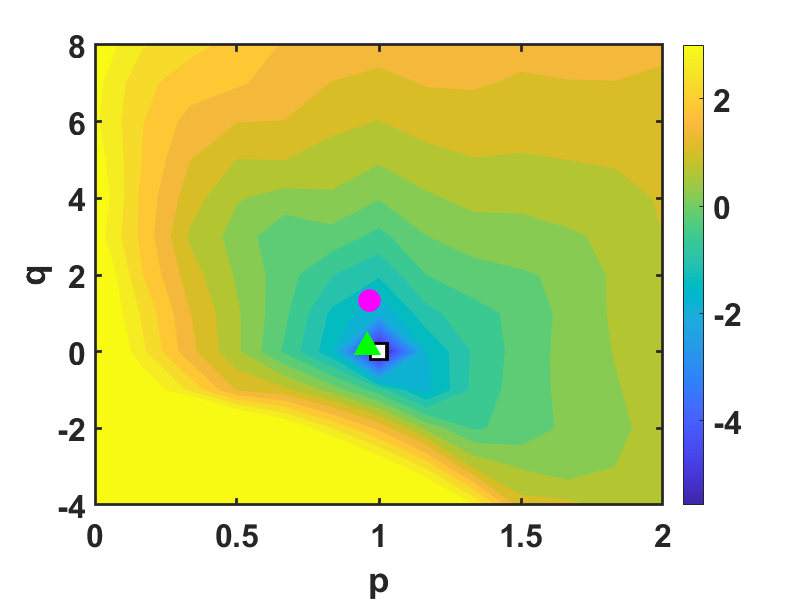}
        \textbf{(a$_\nu$: $\gamma = 0, \nu = 0.001$)}
    \end{minipage}


    \begin{minipage}[c]{0.48\linewidth}
        \centering
        \includegraphics[width=\linewidth,
        trim=0.3cm 0cm 1.3cm 0.8cm, clip]{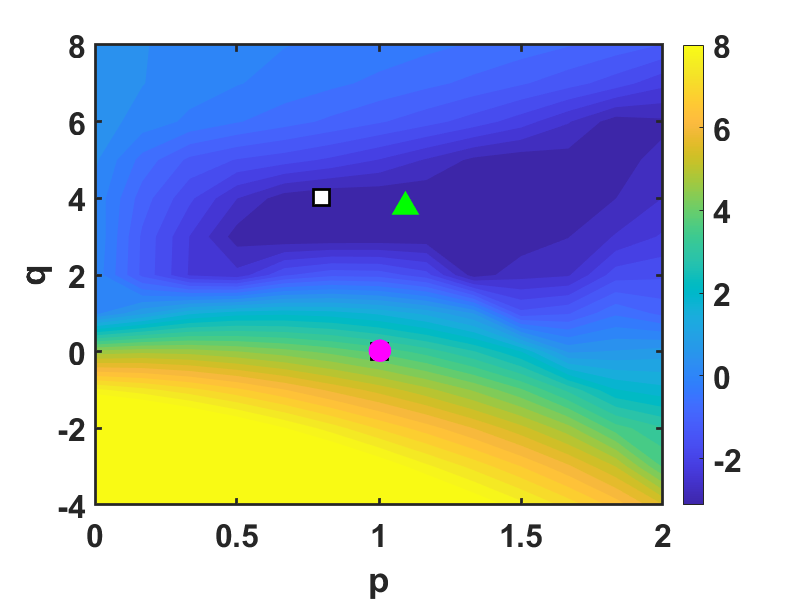}
        \textbf{(b$_\gamma$: $\gamma = 10, \nu = 0$)}
    \end{minipage}%
    \hfill
    \begin{minipage}[c]{0.48\linewidth}
        \centering
        \includegraphics[width=\linewidth,
        trim=0.3cm 0cm 1.3cm 0.8cm, clip]{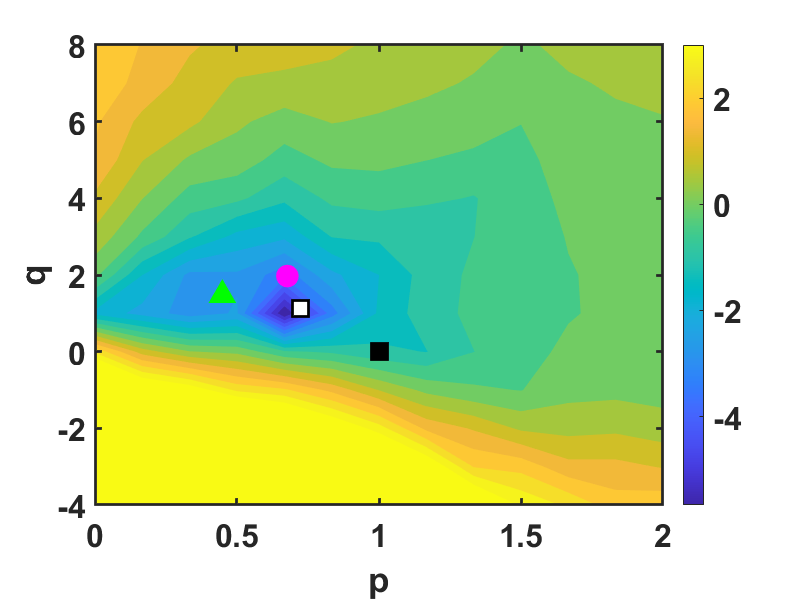}
        \textbf{(b$_\nu$: $\gamma = 0, \nu = 0.05$)}
    \end{minipage}


    \caption{Performance surfaces showing the final SWR error (log-scale) after
    $k = 10$ iterations. Left column: increasing telegrapher damping
    $\gamma > 0$ with $\nu = 0$ (a$_\gamma$: $\gamma = 4$, b$_\gamma$: $\gamma = 10$). Right column: increasing viscoelastic damping
    $\nu > 0$ with $\gamma = 0$ (a$_\nu$: $\nu = 0.001$, b$_\nu$: $\nu = 0.05$).}
    \label{fig:performance_surfaces}
\end{figure}

Figure \ref{fig:performance_surfaces} shows the optimized parameters obtained from the $L_\infty$ (green triangle) and $L_2$ (magenta circle) optimization on the map showing SWR error after 10 iterations.
All optimizations are initialized at $(1/c,0)$ (black square), while the optimal point from the SWR error map is shown as a white square. When $\gamma = 0$ and $\nu > 0$, the computed optimal solutions cluster tightly, and the SWR error map exhibits a pronounced, steep valley around the global minimizer. In contrast, when $\gamma > 0$ and $\nu = 0$, the behavior changes markedly: the error surfaces become much flatter, and the computed minimizers may lie far apart from one another.

\begin{figure}[t]
    \centering


    \begin{minipage}[c]{0.48\linewidth}
        \centering
        \includegraphics[width=\linewidth,
        trim=0.4cm 0cm 0.7cm 0.6cm, clip]{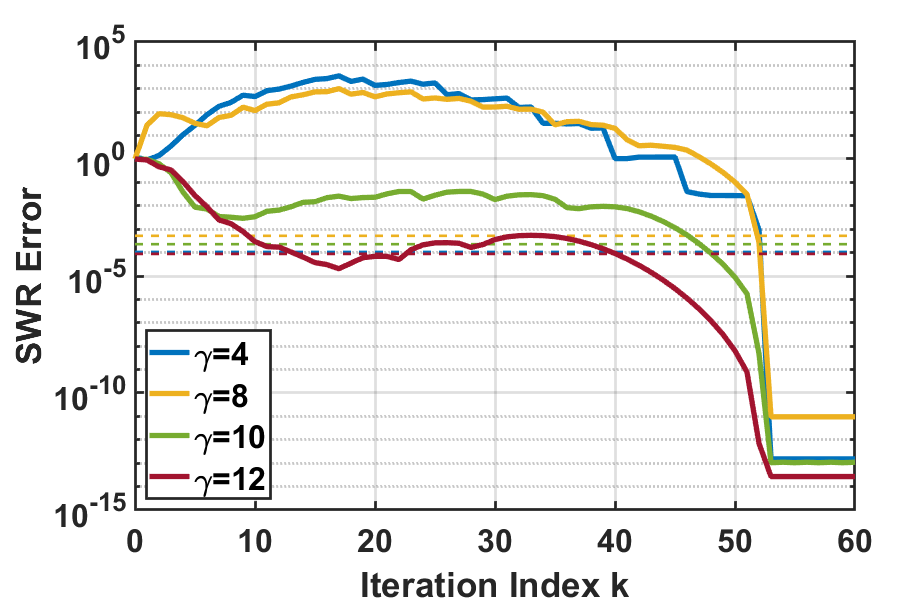}
        \vspace{0.3em}
        \textbf{(a$_\infty$)}
    \end{minipage}%
    \hfill
    \begin{minipage}[c]{0.48\linewidth}
        \centering
        \includegraphics[width=\linewidth,
        trim=0.4cm 0cm 0.7cm 0.6cm, clip]{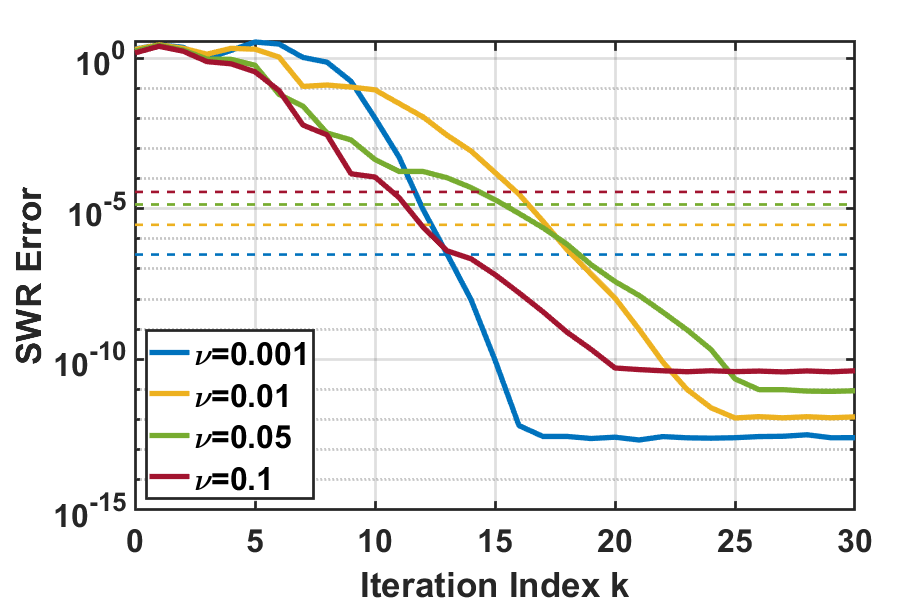}
        \vspace{0.3em}
        \textbf{(b$_\infty$)}
    \end{minipage}



    \begin{minipage}[c]{0.48\linewidth}
        \centering
        \includegraphics[width=\linewidth,
        trim=0.4cm 0cm 0.7cm 0.6cm, clip]{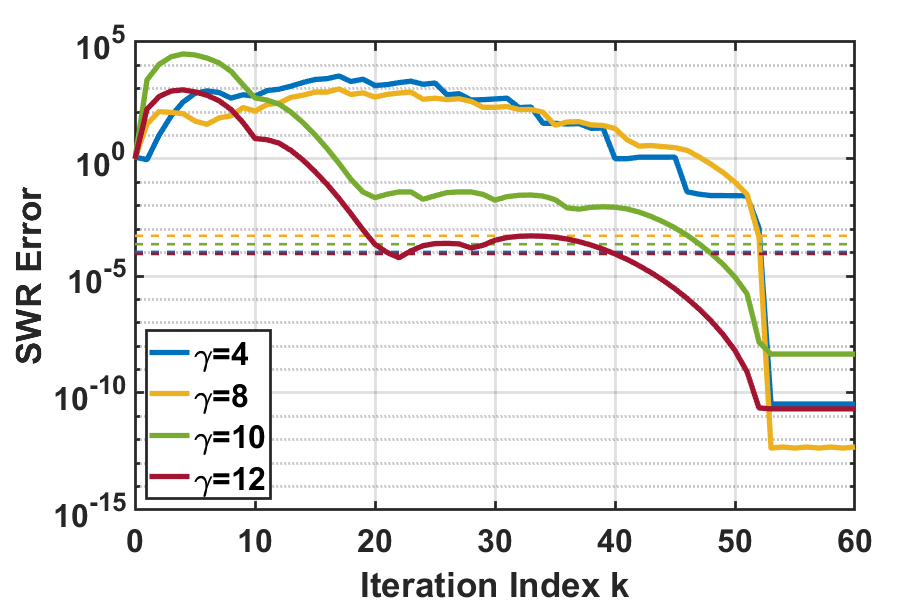}
        \vspace{0.3em}
        \textbf{(a$_2$)}
    \end{minipage}%
    \hfill
    \begin{minipage}[c]{0.48\linewidth}
        \centering
        \includegraphics[width=\linewidth,
        trim=0.4cm 0cm 0.7cm 0.6cm, clip]{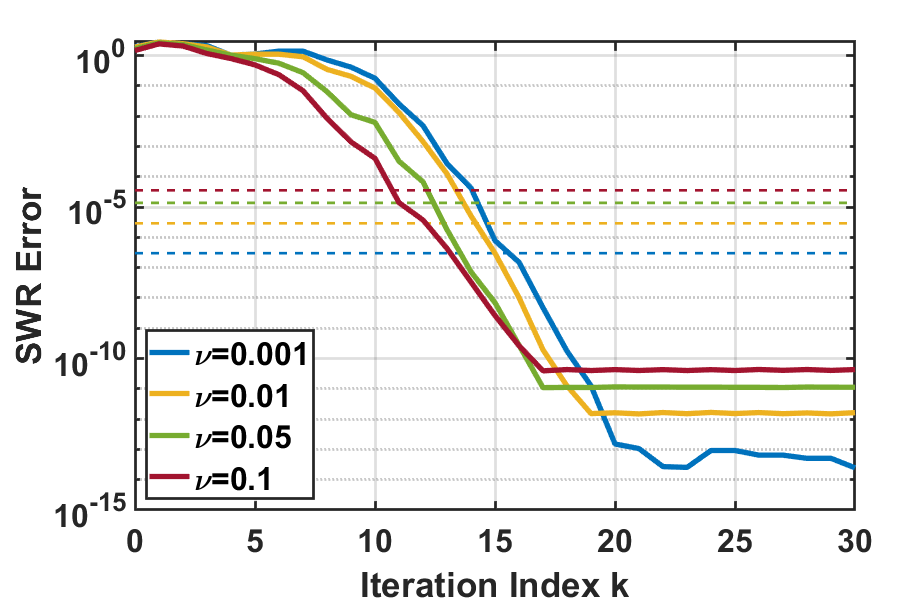}
        \vspace{0.3em}
        \textbf{(b$_2$)}
    \end{minipage}

    \vspace{-0.5em}

    \caption{SWR error curve comparison when varying one damping coefficient at a time. The horizontal axis represents the iteration index k, and the vertical axis represents the relative SWR error on a logarithmic scale.
    First row: $L_\infty$ strategy.
    Second row: $L_2$ strategy.
    (a) Different values of $\gamma$ with $\nu = 0$.
    (b) Different values of $\nu$ with $\gamma = 0$.
    The dashed horizontal line indicates the discrepancy between the full-domain FDTD solution and the analytic 
    ground truth.}
    \label{fig:error_curves_comparison}
\end{figure}

We next plot the error curves in Figure~\ref{fig:error_curves_comparison}.
We begin with the case $\nu = 0$ and vary the telegrapher damping $\gamma$. This corresponds to the first columns of Figure~\ref{fig:error_curves_comparison}.
For $\gamma = 4$ and $\gamma = 8$, Figures~\ref{fig:error_curves_comparison}a show that the SWR error exhibits an initial phase where the error increases; this is followed by a rapid drop to machine precision after about 50 iterations. This behavior is typical of SWR applied to the undamped wave equation, for which it has been proved that the number of iterations required for convergence is at least $cT/\delta$ ~\cite{gander_absorbing_2004,halpern_schwarz_2008}. Although this result has been rigorously established for the undamped case with Dirichlet transmission conditions, the observed behavior here suggests that the telegrapher damping shares a similar convergence mechanism. In fact, the chosen transmission conditions, similarly to Dirichlet ones, do not effectively absorb the incoming errors.
If one increases $\gamma$, the stronger damping mitigates this behavior. In fact, the error curves for $\gamma = 10$ and $\gamma = 12$, using the parameters suggested by the $L_{\infty}$ strategy, decrease at the beginning for about 10 iterations, without the initial growth phase, and then grow a little before, eventually, converging to machine precision after about 50 iterations. The initial decay phase is not present if one considers the parameters suggested by the $L_2$ strategy, which are close to $(1/c,0)$, and do not lie in the optimality region (see Figure \ref{fig:performance_surfaces}).

We conclude that for all telegrapher damping cases, even though the optimized parameters can be close to those minimizing the experimental SWR error (as shown in Fig.~\ref{fig:performance_surfaces}), SWR does not show a clear linear convergence and it is characterized by a convergence mechanism similar to the one of the (undamped) wave equation.

We now turn to the viscoelastic case, where $\nu$ varies while $\gamma = 0$, shown in the second column of Figure~\ref{fig:error_curves_comparison}.  
The error curves in Figures~\ref{fig:error_curves_comparison}b show a systematic improvement in convergence with respect to the telegrapher damping case. The decay is approximately linear, in contrast to the telegrapher case. The parameter optimization is effective. In all cases, the error eventually reaches the machine epsilon. For $\nu = 0.001$, the curve for the $L_{\infty}$ strategy converges faster, while for the three other cases, the $L_2$ strategy provides parameters that exhibit faster convergence.

\begin{figure}[t]
    \centering

    \begin{minipage}[c]{0.48\linewidth}
        \centering
        \vstretch{.9}{\includegraphics[width=\linewidth, trim=1cm 0cm 1cm 0.6cm, clip]{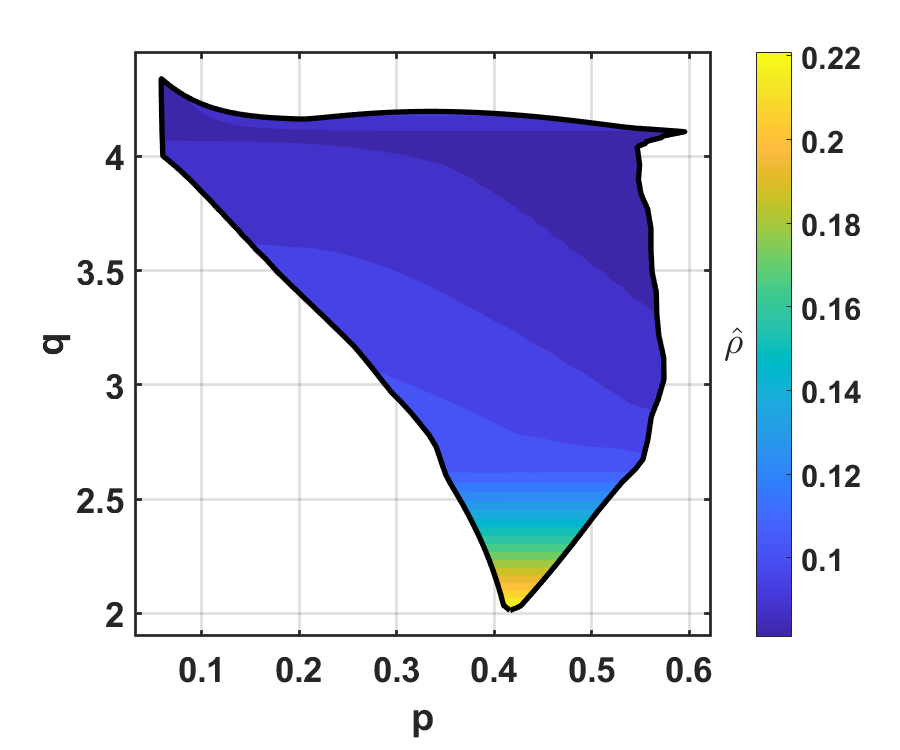}}
        \textbf{(a$_\infty$)}
    \end{minipage}
    \hfill
    \begin{minipage}[c]{0.48\linewidth}
        \centering
        \vstretch{.9}{\includegraphics[width=\linewidth, trim=0.8cm 0cm 1cm 0.6cm, clip]{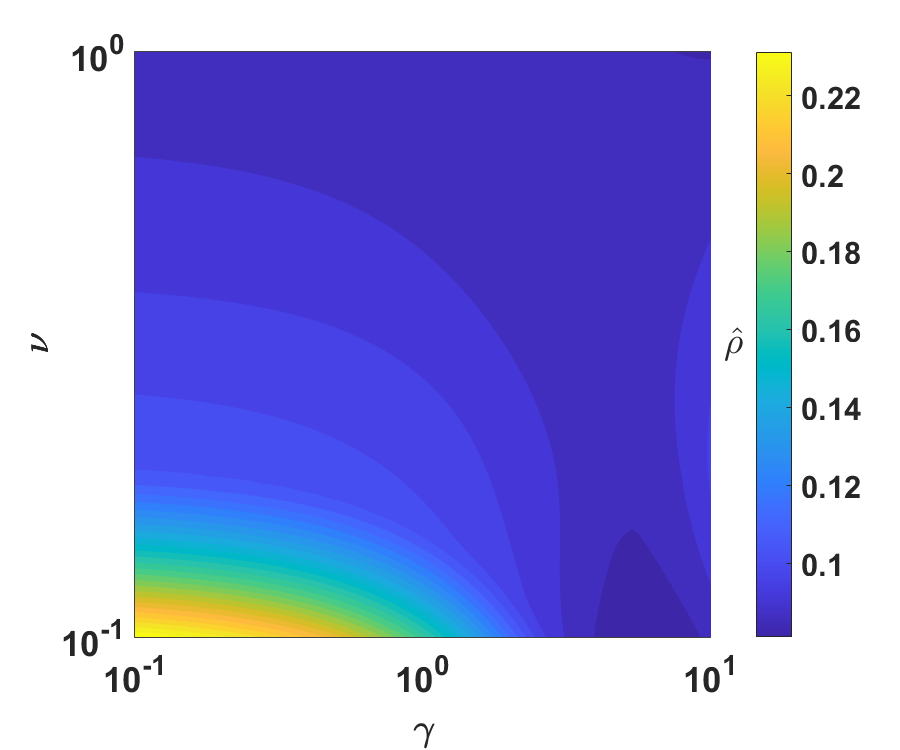}}
        \textbf{(b$_\infty$)}
    \end{minipage}


    
    \begin{minipage}[c]{0.48\linewidth}
        \centering
        \vstretch{.9}{\includegraphics[width=\linewidth, trim=1cm 0cm 1cm 0.6cm, clip]{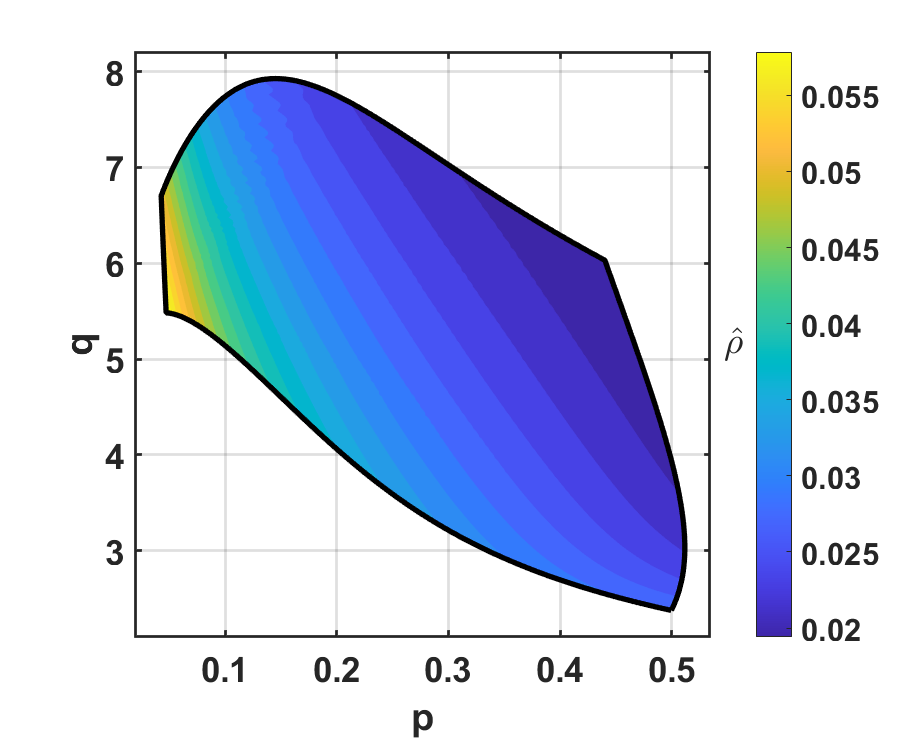}}
        \textbf{(a$_2$)}
    \end{minipage}
    \hfill
    \begin{minipage}[c]{0.48\linewidth}
        \centering
        \vstretch{.9}{\includegraphics[width=\linewidth, trim=0.8cm 0cm 1cm 0.6cm, clip]{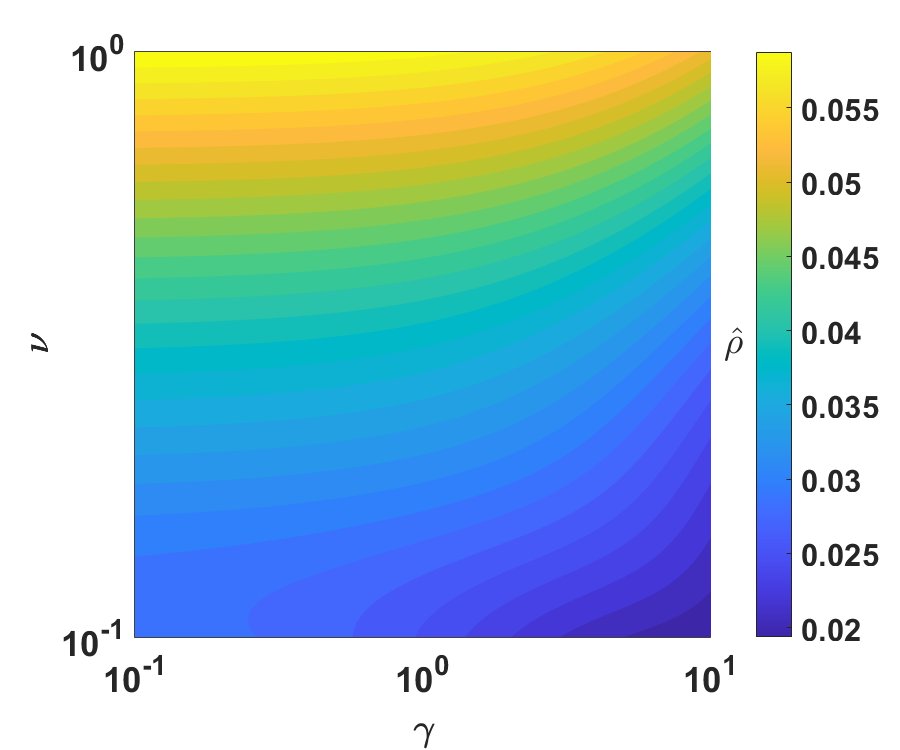}}
        \textbf{(b$_2$)}
    \end{minipage}

    \vspace{-0.5em}

    \caption{Predicted global convergence factors $\hat{\rho}_{\infty}$ (top row) and $\hat{\rho}_2$ (bottom row).
    (a) Contours of $\hat{\rho}$ in the $(p,q)$-plane. 
    (b) Contours of $\hat{\rho}$ in the $(\gamma,\nu)$-plane (log–log axes).}
    \label{fig:comparison_Linf_L2_rho}
\end{figure}

\begin{figure}[t]
    \centering

    \begin{minipage}[c]{0.275\linewidth}
        \centering
        \includegraphics[width=\linewidth]{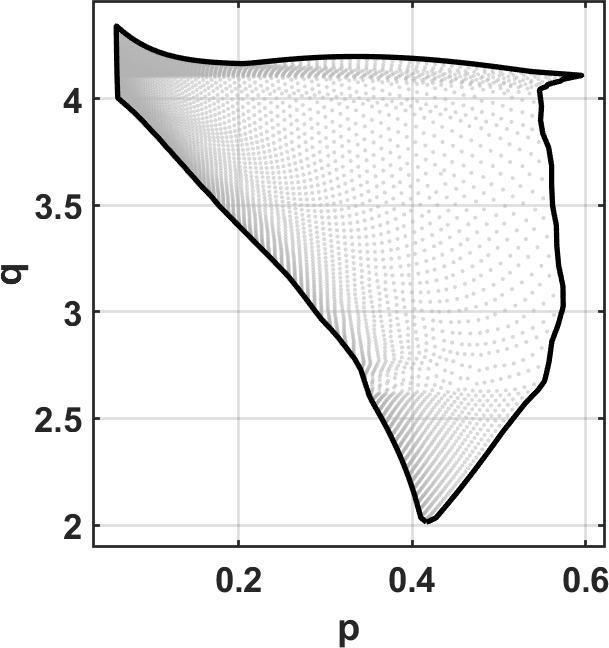}
        \textbf{(c$_\infty$)}
    \end{minipage}%
    \hspace{0.2cm}
    \begin{minipage}[c]{0.34\linewidth}
        \centering
        \includegraphics[width=\linewidth]{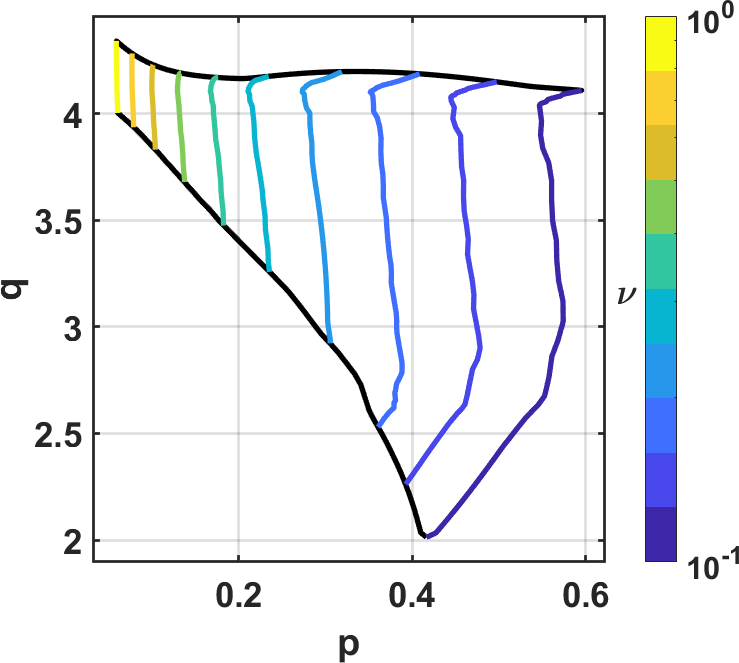}
        \textbf{(d$_\infty$)}
    \end{minipage}%
    \hspace{0.1cm}
    \begin{minipage}[c]{0.34\linewidth}
        \centering
        \includegraphics[width=\linewidth]{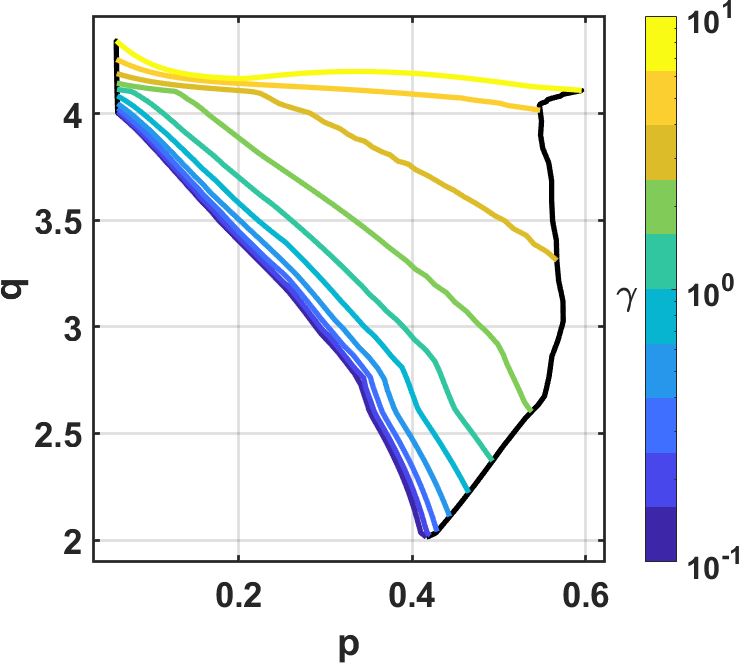}
        \textbf{(e$_\infty$)}
    \end{minipage}


    \begin{minipage}[c]{0.275\linewidth}
        \centering
        \includegraphics[width=\linewidth]{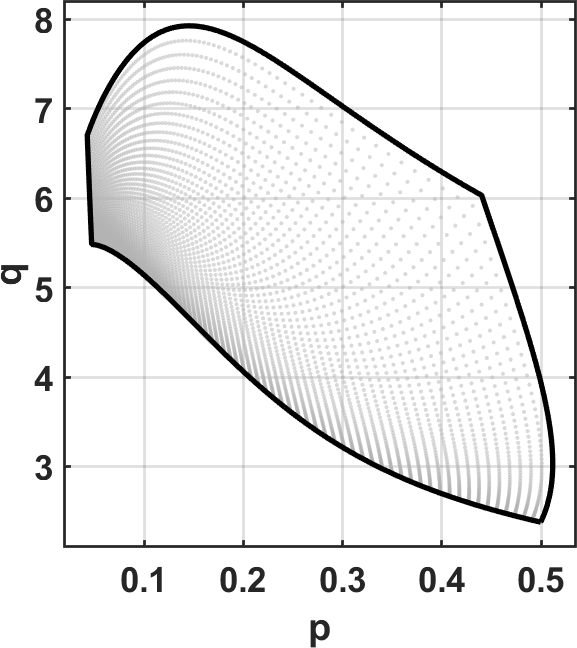}
        \textbf{(c$_2$)}
    \end{minipage}%
    \hspace{0.2cm}
    \begin{minipage}[c]{0.34\linewidth}
        \centering
        \includegraphics[width=\linewidth]{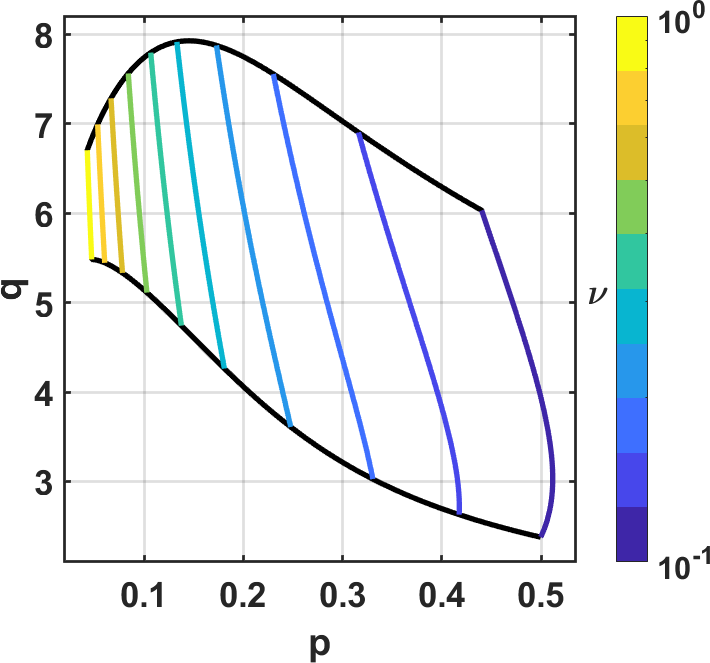}
        \textbf{(d$_2$)}
    \end{minipage}%
    \hspace{0.1cm}
    \begin{minipage}[c]{0.34\linewidth}
        \centering
        \includegraphics[width=\linewidth]{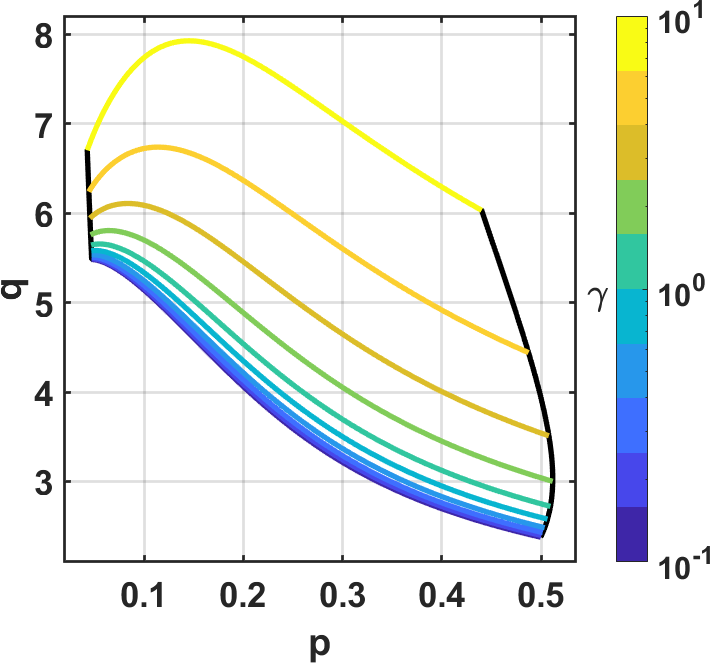}
        \textbf{(e$_2$)}
    \end{minipage}

    \vspace{-0.5em}

    \caption{Comparison between the $L_{\infty}$-optimized (top row) and $L_2$-optimized (bottom row) transmission parameters. 
    (a) Optimal $(p,q)$ points obtained by varying $\gamma$ or $\nu$. 
    (b) Isolines of optimal $(p,q)$ for sweeps at fixed $\nu$. 
    (c) Isolines of optimal $(p,q)$ for sweeps at fixed $\gamma$.}
    \label{fig:comparison_Linf_L2_params}
\end{figure}

In Figures~\ref{fig:comparison_Linf_L2_rho} and \ref{fig:comparison_Linf_L2_params}, we illustrate the dependence of the optimal transmission parameters $(p,q)$ on the damping coefficients $\gamma$ and $\nu$. 
Figure~\ref{fig:comparison_Linf_L2_rho} shows the predicted convergence behavior. Panel~$(a)$ shows the global convergence factor $\hat{\rho}$ in the $(p,q)$–plane, while panel~$(b)$ presents a contour map of $\hat{\rho}$ as a function of $(\gamma,\nu)$, with logarithmic axes. The value of $\hat{\rho}_{\infty}$ decreases when both damping coefficients increase, and variations of $\nu$ produce more pronounced reductions. Instead, $\hat{\rho}_2$ increases as $\nu$ increases.

The first row of Figure~\ref{fig:comparison_Linf_L2_params} (panels $(c_\infty)$ and $(c_2)$) shows how the optimal pairs $(p,q)$ evolve when $\gamma$ varies for several fixed values of $\nu$. The second row (panels $(d_\infty)$ and $(d_2)$) displays isolines obtained by varying $\nu$ at fixed $\gamma$, while the third row (panels $(e_\infty)$ and $(e_2)$) shows isolines produced by varying $\gamma$ at fixed $\nu$. Both optimization strategies agree that the optimal value of $q$ grows as $\gamma$ increases. Under the $L_{\infty}$ criterion, $q$ rises up to approximately $4$ for every tested value of $\nu$, whereas in the $L_2$ case the growth pattern depends on $\nu$, leading to different maximum values of $q$ observed for different damping levels. For both norms, increasing $\nu$ causes the optimal value of $p$ to decrease and converge toward a small value, approximately $p \approx 0.05$.

\section{Conclusions}\label{sec:conclusions}
We presented an optimization framework for SWR transmission conditions applied to the damped wave equation. We showed that the effectiveness of optimized transmission conditions is strictly dependent on the damping model. For the telegrapher equation, the physics of finite wave propagation speed limits the potential for acceleration, resulting in sublinear convergence regardless of parameter tuning. However, for the viscoelastic damped wave equation, the proposed optimization strategies successfully identify transmission parameters that drastically improve performance. The proposed approach offers a computationally efficient alternative to exhaustive search, ensuring robust convergence rates for viscoelastic problems.
Finally, the SWR optimization framework presented here provides a robust basis for accelerating wave simulations in room acoustics and virtual auditory spaces, where real-time performance and computational efficiency are of paramount importance.

\begin{acknowledgement}
The work of G. Ciaramella and I. Mazzieri has been partially supported by the PRIN2022 grant ASTICE - CUP: D53D23005710006.
G. Cicalese, G. Ciaramella, and I. Mazzieri are members of INdAM-GNCS group.
\end{acknowledgement}

\bibliographystyle{spmpsci}
\bibliography{references}

\end{document}